\documentclass{article}

\usepackage{amssymb}
\usepackage{amsfonts}
\usepackage{epsfig}

\newcommand{\Au}{{\bf (A.1)}}
\newcommand{\Ad}{{\bf (A.2)}}
\newcommand{\At}{{\bf (A.3)}}
\newcommand{\Aq}{{\bf (A.4)}}
\newcommand{\tu}{T_n^{(1)}}
\newcommand{\tr}{T_n^{(2,\rho)}}

\newcommand{\trz}{T_n^{(2,0)}}
\newcommand{\tru}{T_n^{(2,\rho,1)}}
\newcommand{\trd}{T_n^{(2,\rho,2)}}

\newcommand{\tuu}{T_n^{(1,1)}}
\newcommand{\tud}{T_n^{(1,2)}}

\newcommand{\maxE}{\|\varepsilon\|_{n,\infty}}

\newcommand{\ain}{\alpha_{i,n}}
\newcommand{\bin}{\alpha_{k-i,n}}
\newcommand{\ein}{\varepsilon_{i,n}}
\newcommand{\fin}{\varepsilon_{k-i,n}}
\newcommand{\unchap}{X_{n-k_n+1,n}}
\newcommand{\somi}{\sum_{i=1}^{k-1}}
\newcommand{\somii}{\sum_{i=1}^{k_n-1}}
\newcommand{\invH}{H^{\leftarrow}}

\newcommand{\tod}{\stackrel{d}{\to}}

\newcommand{\egl}{\stackrel{d}{=}}
\newcommand{\ld}[1] {\log\log\left({#1}\right)}
\newcommand{\CQFD}
{%
\mbox{}%
\nolinebreak%
\hfill%
\rule{2mm}{2mm}%
\medbreak%
\par%
}
\newcommand{\proof}{{\bf Proof~: }}

\newtheorem{Theo}{Theorem}
\newtheorem{Coro}{Corollary}
\newtheorem{Lem}{Lemma}

\begin{document}

\title{Estimation of the Weibull tail-coefficient with linear combination
of upper order statistics}
\author{Laurent Gardes \& St\'ephane Girard}
\date{INRIA Rh\^one-Alpes, projet Mistis,
Inovall\'ee, 655, av. de l'Europe, Montbonnot, 38334 Saint-Ismier cedex,
France.
{\tt \{Laurent.Gardes,Stephane.Girard\}@inrialpes.fr}}

\maketitle

\begin{abstract}
We present a new family of estimators of the Weibull tail-coefficient. The Weibull
tail-coefficient is defined as the regular variation coefficient of
the inverse failure rate function. Our estimators are based on a linear combination of log-spacings of the upper order statistics. Their asymptotic normality is established and illustrated for two
particular cases of estimators in this family.  Their finite sample performances are presented on a simulation study.\\\\

\noindent{\bf Keywords}: Weibull tail-coefficient, extreme-values, order statistics, regular variations.
\end{abstract}

\section{Introduction}
\label{intro}

Weibull tail-distributions encompass a variety of light tailed distributions,
{\it i.e.} distributions in the Gumbel maximum domain of attraction,
see~\cite{Galambos} for further details. Weibull tail-distributions include
for instance Weibulls, Gaussians, gammas and logistic.
The purpose of this paper is to study the estimation of a tail parameter
associated with these distributions.
More precisely, a cumulative distribution function $F$ has a Weibull tail
if its logarithmic tail satisfies the following property:
There exists $\theta>0$
such that for all $\lambda>0$,
\begin{equation}
\label{A1bis}
\lim_{t\to\infty} \frac{\log(1-F(\lambda t))}{\log(1-F(t))}=\lambda^{1/\theta}.
\end{equation} 
The parameter of interest $\theta$ is called the Weibull tail-coefficient.
Such distributions are of great use to model large claims in non-life
insurance~\cite{BeirTeug}. 
In the particular case where
$
{\log(1-F(\lambda t))}/{\log(1-F(t))}=\lambda^{1/\theta}
$
for all $t>0$ and $\lambda>0$,
estimating $\theta$ reduces to
estimating the shape parameter of a Weibull distribution.
In this context, simple and efficient methods exist, see for
instance~\cite{BAIN}, Chapter~4 for a review on this topic.
Otherwise,
dedicated estimation methods have been proposed
since the relevant information on the Weibull tail-coefficient
is only contained in the extreme upper part of the sample.
A first direction was investigated in~\cite{Berred} where an
estimator based on the record values is proposed.
Another family of approaches~\cite{BRO,BeirBro,BTV,thetaWT}
consists of using the $k_n$ upper order statistics where
$(k_n)$ is an intermediate sequence of integers {\it i.e.}
such that
\begin{equation}
\label{intermed}
\lim_{n\to\infty} k_n=\infty \mbox{ and } \lim_{n\to\infty} k_n/n=0.
\end{equation}
Note that, since $\theta$ is defined only by an asymptotic behavior
of the tail, the estimator should use the only extreme-values of the
sample and thus the second part of~(\ref{intermed}) is required.
The estimators considered here belong to this approach.
Let $(X_i)_{1\leq i\leq n}$ be a sequence of independent and identically
distributed random variables with cumulative distribution function~$F$.
Denoting by $X_{1,n}\leq\dots\leq X_{n,n}$ the corresponding order statistics,
our family of estimators is
\begin{equation}
\label{deftheta}
\hat{\theta}_n(\alpha) =
{\displaystyle\somii \ain(\log(X_{n-i+1,n})  - \log(\unchap)) }
\left/{\displaystyle\somii \alpha_{i,n}(\ld{n/i} - \ld{n/k_n})} 
\right.  
\end{equation}
with weights $\ain=W(i/k_n) + \ein$ defined from
$W$ a smooth score function
and $(\ein)_{1\leq i\leq k_n-1}$ a non-random sequence.
We refer to~\cite{CDM,Viharos} for similar works in the
context of the estimation of the extreme-value index.

In Section~\ref{main} we state the asymptotic normality of these estimators.
In Section~\ref{comp}, we provide two examples of weights.
The first one leads to the estimator of $\theta$ proposed
by Beirlant {\it et al.}~\cite{BTV}. The second one gives rise to a new estimator
for Weibull tail-distributions. The behavior of these two estimators
is investigated on finite sample situations in~\cite{Nous}.
Finally, proofs are given in Section~\ref{preuves}.

\section{Asymptotic normality}
\label{main}

Consider the failure rate $H=-\log(1-F)$. Writing $\invH$ its
generalized inverse $\invH(t)=\inf\{x,\; H(x)\geq t\}$,
assumption~(\ref{A1bis}) is equivalent to:
\begin{description}
\item [(A.1)] $\invH(t)=t^\theta\ell(t)$,
\end{description}
\noindent where $\ell$ is a slowly varying function {\it i.e.} such that
$\ell(\lambda t)/\ell(t)\to 1$ as $t\to\infty$ for all $\lambda>0$.
The inverse failure rate function $\invH$ is said to be
regularly varying at infinity with index~$\theta$
and this property is denoted by $\invH\in{\mathcal R}_\theta$.
We refer to~\cite{BING} for more information on regular variation theory.  
As a comparison, Pareto type distributions satisfy
$(1/(1-F))^\leftarrow \in{\mathcal R}_\gamma$, and
$\gamma>0$ is the so-called extreme-value index.
As often in extreme-value theory, \Au~is not sufficient
to prove a central limit theorem for $\hat\theta_n(\alpha)$.
It needs to be strengthened with a second order condition on $\ell$,
namely that there exist $\rho\leq 0$ and a function $b$ with
limit 0 at infinity such that 
\begin{description}
\item [(A.2)] 
$
 \log\left({\ell(\lambda t)}/{\ell(t)}\right) \sim b(t)\int_1^\lambda u^{\rho-1} du , 
$
\end{description}
uniformly locally on $\lambda>1$ and as $t\to\infty$.
The second order parameter $\rho\leq 0$ tunes the rate of convergence
of $\ell(\lambda t)/\ell(t)$ to 1. The closer $\rho$ is to 0, the slower
is the convergence. Condition \Ad~is the cornerstone in all proofs of
asymptotic normality for extreme-value estimators.
It is used in~\cite{Hill,Hausler,BEIR} to prove the asymptotic normality of
estimators of the extreme-value index~$\gamma$. 
Table~\ref{tabex} shows that many distributions satisfy \Au~and \Ad.
Among them, Extended Weibull distributions, introduced in~\cite{KLUP},
encompass gamma, Gaussian and Benktander~II distributions.
We refer to~\cite{EMBR}, Table~3.4.4, for the derivation of $b(x)$ and $\rho$
in each case.
Other examples are the Weibull, logistic and extreme-value 
(with shape parameter $\gamma=0$)
distributions.

Throughout the paper, we write Id for the identity function. In particular,
if $f$ is a function and $p$ a real number, the inequality $f\leq \mbox{Id}^p$
means $f(t)\leq t^p$ for any $p$ where it is defined.
For general L-estimators, conditions on the weights are
required to obtain a central limit theorem
(see for instance~\cite{Mason2}). 
Our assumptions are the following:
\begin{description}
\item [(A.3)] $W$ is defined and continuously differentiable on the open unit interval,
\item [(A.4)] There exist $M>0$, $0\leq q< 1/2$ and $p<1$
such that $|W|\leq M \mbox{Id}^{-q}$
and $|W'|\leq M \mbox{Id}^{-p-q}$ on the open unit interval. 
\end{description}
Similar conditions have been introduced in the
context of the estimation of the extreme-value index~\cite{CDM,Viharos}.
To write the limiting variance of $\hat{\theta}_n(\alpha)$, we introduce two
quantities:
\begin{eqnarray*}
\mu(W)&=&\int_0^1 W(x)\log(1/x)dx, \\
\sigma^2(W)&=&\int_0^1\int_0^1  W(x)W(y)\frac{\min(x,y)-xy}{xy}dxdy.
\end{eqnarray*}
We also define $\maxE=\displaystyle\max_{i=1,\dots,k_n-1}|\ein|$.  
We are now in position to state our main result.
Its proof is postponed to Section~\ref{preuves}. 
\begin{Theo}
\label{normalite}
Suppose \Au--\Aq~hold. 
If $(k_n)$ is any intermediate sequence such that
\begin{equation}
\label{condkn}
  k_n^{1/2}b(\log(n))\to\lambda \mbox{ and }
k_n^{1/2}\max\{1/\log(n), \maxE\}\to 0,
\end{equation}
then
$$
k_n^{1/2}(\hat{\theta}_n(\alpha)-\theta) \tod {\mathcal N}(\lambda,\theta^2\sigma^2(W)/\mu^2(W)).
$$
\end{Theo}
Clearly, the bias of the estimator is driven by the function $b$.  
This bias term asymptotically vanishes if $\lambda=0$.
Some applications of this result are 
given in the next section, Corollary~\ref{coro1} and 
Corollary~\ref{normalitezipf}.
The importance of the bias term is also illustrated on
finite sample situations.  
Finally, note that condition (\ref{condkn}) implies $k_n/n\to 0$.

\section{Comparison of two estimators}
\label{comp}

First, we show in Paragraph~\ref{casgirard}, that our family of estimators~(\ref{deftheta}) encompasses
the Hill type estimator $\hat{\theta}_n^{H}$ proposed in~\cite{BTV}. Moreover,
it will appear in Corollary~\ref{coro1} that the asymptotic normality
of $\hat{\theta}_n^{H}$ stated in~\cite{thetaWT}, Theorem~2 is a
consequence of our main result Theorem~1.
Second, in Paragraph~\ref{caszipf}, we use our framework to 
exhibit a new estimator of the Weibull tail-coefficient and to establish
its asymptotic normality in Corollary~\ref{normalitezipf}.  
In the third paragraph, we show that the new estimator performs
as well as the Hill one.

\subsection{Hill type estimator}
\label{casgirard}

Beirlant {\it et al.}~\cite{thetaWT} propose the following estimator of the Weibull tail-coefficient:
\begin{equation}
\label{defthetasteph}
\hat{\theta}_n^{H} =
{\displaystyle\somii (\log(X_{n-i+1,n})  - \log(\unchap)) } 
\left/{\displaystyle\somii (\ld{n/i} - \ld{n/k_n})}\right..
\end{equation}
Clearly, $\hat{\theta}_n^{H}$ is a particular case of $\hat\theta_n(\alpha)$ 
with $W(x) = 1$ for all $x\in[0,1]$ and $\ein = 0$ for all $i=1,\ldots,k_n$. 
The asymptotic normality of $\hat{\theta}_n^{H}$, established
in Theorem~2 of~\cite{thetaWT}, can be obtained as a
consequence of Theorem~\ref{normalite}:
\begin{Coro}
\label{coro1}
Suppose \Au~and \Ad~hold. If $(k_n)$ is an
 intermediate sequence such that
$
k_n^{1/2}\max\{b(\log(n)), 1/\log(n)\}\to 0,
$
then $
k_n^{1/2} (\hat{\theta}_n^{H}-\theta) \tod {\mathcal N}(0,\theta^2).
$
\end{Coro}

\subsection{Zipf estimator}
\label{caszipf}

We propose a new estimator of the Weibull tail-coefficient
based on a quantile plot adapted to our situation.
It consists of drawing the pairs $(\ld{n/i},\log{(X_{n-i+1,n})})$
for $i=1,\dots,n-1$. The resulting graph
should be approximatively linear (with slope $\theta$),
at least for the large values of $i$.
Thus, we introduce $\hat{\theta}_n^{Z}$ the least square estimator of $\theta$ based on
the $k_n$ largest observations:
\begin{equation}
\label{defthetazipf}
\hat{\theta}_n^{Z} =
{\displaystyle\somii  ( \ld{n/i}-\zeta_n)\log(X_{n-i+1,n}) } 
\left/{\displaystyle\somii ( \ld{n/i}-\zeta_n)  \ld{n/i} } \right.,
\end{equation}
where
$$
\zeta_n = \frac{1}{k_n-1} \somii \ld{n/i}.
$$
This estimator is similar to the Zipf estimator for the extreme-value index proposed by
Kratz and Resnick~\cite{krares96} and Schultze and Steinebach~\cite{schste96}. 
We prove in Section~\ref{preuves}
that $\hat{\theta}_n^{Z}$ belongs to family~(\ref{deftheta})
and thus apply Theorem~\ref{normalite}
to obtain its asymptotic normality:
\begin{Coro}
\label{normalitezipf}
Suppose \Au~and \Ad~hold. If $(k_n)$ is an
 intermediate sequence such that
$k_n^{1/2}\max\{b(\log(n)), \log^2(k_n)/\log(n)\}\to 0,
$
then $
k_n^{1/2} (\hat{\theta}_n^{Z}-\theta) \tod {\mathcal N}(0,2\theta^2).
$
\end{Coro}

\section{Proofs}
\label{preuves}

Throughout this section, we assume that $(k_n)$ is an intermediate sequence
and, for the sake of simplicity, we note $k$ for $k_n$. 
Let us also introduce $K_\rho(\lambda)=\int_1^\lambda u^{\rho-1} du$ for 
$\lambda\geq 1$
and $J(x)=W(1-x)$ for $x\in(0,1)$.  
The following notations will prove useful: 
 $E_{n-k+1,n}$ is the $(n-k+1)$th order statistics associated 
to  $n$ independent standard exponential variables and
 $(F_{i,k-1})_{1\leq i\leq k-1}$ are order statistics,
independent from $E_{n-k+1,n}$, 
generated by $k-1$ independent standard exponential variables.
The next lemma presents an expansion of $\hat\theta_n(\alpha)$.  
\begin{Lem}
\label{lemdecomp}
Under~\Au~and~\Ad, $\hat\theta_n(\alpha)$ has the same distribution as
$$
 \frac{\theta\trz+(1+o_P(1))b(E_{n-k+1,n}) \tr
}{\tu} ,
$$
where we have defined
\begin{eqnarray*}
\tu&=&\frac{1}{k-1} \somi \ain (\ld{n/i} - \ld{n/k}) \mbox{ and }\\
\tr&=& \frac{1}{k-1}\somi \bin K_\rho\left(1+\frac{F_{i,k-1}}{E_{n-k+1,n}}\right),
\; \rho\leq 0.
\end{eqnarray*}
\end{Lem}
\proof
Using the quantile transform, the order statistics $(X_{i,n})_{1\leq i\leq n}$
have the same distribution as $(\invH(E_{i,n}))_{1\leq i\leq n}$.
Thus, \Ad~yields that the numerator of $\hat\theta_n(\alpha)$ in~(\ref{deftheta})
has the same distribution as
$$
 \theta\frac{1}{k-1}\somi \ain \log\left(\frac{E_{n-i+1,n}}{E_{n-k+1,n}}\right)
+ (1+o_P(1)) b(E_{n-k+1,n}) \frac{1}{k-1}\somi \ain K_\rho\left(\frac{E_{n-i+1,n}}{E_{n-k+1,n}}\right).
$$
The R\'enyi representation  asserts that
$ ({E_{n-i+1,n}}/{E_{n-k+1,n}})_{1\leq i\leq k-1}$
has the same distribution as
$( 1+ {F_{k-i,k-1}}/{E_{n-k+1,n}})_{1\leq i\leq k-1}$,
see~\cite{order}, p.~72.
Therefore, the numerator of $\hat\theta_n(\alpha)$ has the same distribution as
\begin{eqnarray*}
&& \theta\frac{1}{k-1} \somi \ain \log\left(1+\frac{F_{k-i,k-1}}{E_{n-k+1,n}}\right)\\
&+& (1+o_P(1)) b(E_{n-k+1,n}) \frac{1}{k-1}\somi \ain K_\rho\left(1+\frac{F_{k-i,k-1}}{E_{n-k+1,n}}\right).
\end{eqnarray*}
Changing $i$ to $k-i$ in the above formula and remarking that
$K_0$ is the logarithm function conclude the proof.
\CQFD
\noindent The following lemma provides an expansion of
$$
\tau_n = \frac{1}{k-1} \somi (\ld{n/i}-\ld{n/k}),
$$
which frequently appears in the proofs.
\begin{Lem}
\label{lemtaun}
The following expansion holds:
$$
\tau_n = \frac{1}{\log(n/k)} \left\{ 1 + O\left(\frac{\log(k)}{k}\right) + O\left ( \frac{1}{\log(n/k)} 
\right) \right\} .
$$
\end{Lem}
\proof
We write $\tau_n$ as the sum 
$$
\frac{1}{\log(n/k)} \frac{1}{k-1} \somi \log(k/i) 
+
\frac{1}{k-1} \somi \left \{ \log \left ( 1+\frac{\log(k/i)}{\log(n/k)} \right ) -\frac{\log(k/i)}{\log(n/k)} \right \}.
$$
Since
$$
 \frac{1}{k-1} \somi \log(i/k)=\frac{1}{k-1} \log \left( \frac{k!}{k^k}\right),
$$
Stirling's formula shows that the first term is
$$
\frac{1}{\log(n/k)}\left(1 +  O\left(\frac{\log(k)}{2k}\right)\right).
$$
The inequality $-x^2/2\leq \log(1+x)-x \leq 0$, valid for nonnegative $x$
shows that the second term is of order at most
$$
\frac{1}{\log^2(n/k)} \frac{1}{k-1} \somi \log^2(i/k) = O \left ( \frac{1}{\log^2(n/k)} \right ),
$$
since the above Riemann sum converges to 2
as $k \to \infty$. The result follows.
\CQFD
\noindent
The next lemmas are dedicated to the study of the
different terms appearing in Lemma~\ref{lemdecomp}.
First, we focus on the non-random term $\tu$.  
\begin{Lem}
\label{lemTun}
Under~\Au--\Aq, the following expansion hold: 
$$
\tu=
\frac{\mu(W)}{\log (n/k)}
\left\{
1+O\left(\log(k)k^{q-1}\right)
+O\left(\frac{1}{\log{(n/k)}}\right)
+O\left(\maxE\right)
\right\}.
$$
\end{Lem}
\proof
Clearly, $\tu$ can be rewritten as the sum
$$
\frac{1}{k-1} \somi \ein \log\left(1+\frac{\log(k/i)}{\log(n/k)}\right)
+\frac{1}{k-1} \somi W(i/k) \log\left(1+\frac{\log(k/i)}{\log(n/k)}\right).
$$
The absolute value of the first term is less than 
$\maxE \tau_n$ which is $O({\maxE}/{\log(n/k)})$,
by Lemma~\ref{lemtaun}.
The second term can be expanded as
\begin{eqnarray*}
\frac{1}{\log(n/k)}\frac{1}{k-1} \somi W(i/k) \log(k/i) 
+ \frac{1}{k-1} \somi W(i/k) \left\{
\log\left(1+\frac{\log(k/i)}{\log(n/k)}\right) -
 \frac{\log(k/i)}{\log(n/k)}\right\}&&\\
=: \frac{\tuu}{\log(n/k)} + \tud.&&
\end{eqnarray*}
For $x\in(0,1)$, define $H(x)=W(x) \log(1/x)$.
The Riemann sum $\tuu$ can be compared to $\mu(W)$ by:
\begin{equation}
\label{firstterm}
| \tuu - \mu(W) | \leq 
\frac{1}{2k^2} \somi \sup_{i/k \leq x \leq (i+1)/k} | H'(x) | + \int_0^{1/k} | H(x) | dx + O(1/k). 
\end{equation}
Assumption~\Aq~implies that there exists a positive $M'$ such that
$|H'|\leq M' \mbox{Id}^{-q-1}$ on the open unit interval, and thus
the first term of~(\ref{firstterm}) is bounded above by
$$
\frac{M'}{2k} \left ( \int_{1/k}^1 t^{-q-1}dt +  k^q \right) = \left \{ \begin{array}{l
l}
O\left(k^{q-1}\right) & {\rm{if}} \ q \neq 0, \\
O\left(k^{-1}\log(k)\right) & {\rm{otherwise}}.
\end{array} \right .
$$
Assumption~\Aq~also yields $|H|\leq M\mbox{Id}^{-q}\log(1/x)$ on the open unit interval
and thus the second term in~(\ref{firstterm}) is 
$O\left(k^{q-1}\log (k)\right)$.
It follows that
\begin{equation}
\label{termtuu}
\tuu = \mu(W)+O(k^{q-1}\log(k)).
\end{equation}
Besides, the well-known inequality $|\log (1+x) -x | \leq x^2/2$, valid for all nonnegative
$x$ together with~\Aq~show that $|\tud|$ is bounded by
\begin{equation}
\label{termtud}
|\tud| 
 \leq \frac{M}{2\log^2(n/k)} \frac{1}{k-1}\somi (i/k)^{-q}\log^2(k/i)
 = O\left(\frac{1}{\log^2(n/k)}\right),
\end{equation}
 since the above Riemann sum converges to a finite integral.
Collecting (\ref{termtuu}) and (\ref{termtud}) gives
the result.  \CQFD
\noindent Second, we focus on the random term $\tr$.  
\begin{Lem}
\label{lemTtrois}
Let
$\xi$ be standard Gaussian random variable. 
Under~\Au--\Aq, the following expansion hold for all non-positive $\rho$:
$$
\tr\egl 
\frac{\mu(W)}{E_{n-k+1,n}}
\left\{
1+ 
\frac{\sigma(W)}{\mu(W)} k^{-1/2}\xi(1+o_P(1))
+O_P\left(\frac{1}{\log{(n/k)}}\right)
+O_P\left(\maxE\right)
\right\}.
$$
\end{Lem}
\proof
Note that $\tr$ can be written as the sum
\begin{equation}
\label{eqdebut}
\frac{1}{k-1} \somi \fin K_\rho\left(1+\frac{F_{i,k-1}}{E_{n-k+1,n}}\right)
+\frac{1}{k-1} \somi J(i/k) K_\rho\left(1+\frac{F_{i,k-1}}{E_{n-k+1,n}}\right). \end{equation} 
Since $0\leq K_\rho(1+x)\leq x$ for all nonnegative $x$, the absolute value
of the first term is bounded by
$$
\maxE \frac{1}{k-1} \somi K_\rho\left(1+\frac{F_{i,k-1}}{E_{n-k+1,n}}\right)
\leq \maxE \frac{1}{k-1}  \somi  \frac{F_{i,k-1}}{E_{n-k+1,n}},
$$
which has the same distribution as 
$$
\frac{\maxE}{E_{n-k+1,n}}  \frac{1}{k-1} \somi F_i = 
\frac{1}{E_{n-k+1,n}} O_P\left(\maxE\right),
$$
from the law of large numbers.   
The second term of~(\ref{eqdebut}) can be expanded as
\begin{eqnarray*}
\frac{1}{E_{n-k+1,n}}\frac{1}{k-1} \somi J(i/k) F_{i,k-1}
+ \frac{1}{k-1} \somi J(i/k) \left\{
K_\rho\left(1+\frac{F_{i,k-1}}{E_{n-k+1,n}}\right) -
 \frac{F_{i,k-1}}{E_{n-k+1,n}}\right\}&&\\
=: \frac{\tru}{E_{n-k+1,n}} + \trd.&&
\end{eqnarray*}
Now, \At~and~\Aq~imply that the L-statistics
$\tru$ satisfies the conditions of~\cite{Mason2}  and thus is asymptotically
Gaussian.   More precisely, we have
\begin{equation}
\label{eqtru}
\tru\egl\mu(W)+ \sigma(W) k^{-1/2}\xi(1+o_P(1)).
\end{equation}
The upper bound on $\trd$ is obtained by remarking that
$
| K_\rho(1+x) -x | \leq (1-\rho) x^2/2
$
for all nonnegative $x$.
It follows that $\trd$ is bounded above by
$$
\frac{1}{k-1} \somi |J(i/k)| \left|
K_\rho\left(1+\frac{F_{i,k-1}}{E_{n-k+1,n}}\right) -
 \frac{F_{i,k-1}}{E_{n-k+1,n}}\right|
 \leq \frac{1-\rho}{2E_{n-k+1,n}^2}  \frac{1}{k-1}  \somi |J(i/k)| F_{i,k-1}^2.
$$
Now, $E_{n-k+1,n}$ is equivalent to $\log{(n/k)}$ in probability and
$$
\frac{1}{k-1} \somi |J(i/k)| F_{i,k-1}^2 =O_P(1),
$$
from the results of~\cite{Mason2} on L-statistics.
Thus
\begin{equation}
\label{eqtrd}
\trd=\frac{1}{E_{n-k+1,n}}O_P\left(\frac{1}{\log{(n/k)}}\right),
\end{equation}
and then collecting (\ref{eqtru}) and (\ref{eqtrd}), the second term 
of~(\ref{eqdebut}) is
$$
\frac{1}{E_{n-k+1,n}} \left ( \mu(W)+\sigma(W)k^{-1/2}\xi(1+o_P(1))+O_P\left(\frac{1}{\log{(n/k)}} \right ) \right ),
$$
and the result follows.  \CQFD

\noindent We are now in position to prove Theorem~\ref{normalite} and Corollary~\ref{normalitezipf}.
 
\paragraph{\bf Proof of Theorem \ref{normalite}.}

From Lemma~\ref{lemdecomp}, $k^{1/2}(\hat{\theta}_n(\alpha)-\theta)$ has the
same distribution as
$$
\theta k^{1/2}\left(\frac{\trz}{\tu} -1 \right) + 
k^{1/2} b(E_{n-k+1,n})\frac{\tr}{\tu}(1+o_P(1)).
$$
Now, $E_{n-k+1,n}$ is equivalent to $\log(n/k)$ in probability
which is also equivalent to $\log(n)$, see Lemma~5.1 in~\cite{Revstat}.
Since $|b|$ is regularly varying (see \cite{GEL}),
 $b(E_{n-k+1,n})$ is equivalent to $b(\log(n))$
in probability. As a consequence, 
$k^{1/2}(\hat{\theta}_n(\alpha)-\theta)$ has the
same distribution as
\begin{equation}
\label{eqtmp1}
\theta k^{1/2}\left(\frac{\trz}{\tu} -1 \right) +
k^{1/2} b(\log(n))\frac{\tr}{\tu}(1+o_P(1)).
\end{equation}
Let us consider the first term of this sum.
Lemma~\ref{lemTun}, Lemma~\ref{lemTtrois} and condition~(\ref{condkn}) entail that, for all non-positive $\rho$,
the ratio $\tr/\tu$ has the same distribution as
$$
\frac{\log(n/k)}{E_{n-k+1,n}}\left\{
 1   +  \frac{\sigma(W)}{\mu(W)} k^{-1/2} \xi(1+o_P(1)) \right\}.
$$
Now, Lemma~1 in~\cite{thetaWT} asserts a central limit theorem for order
statistics of an exponential sample, and thus
$$
\frac{\log(n/k)}{E_{n-k+1,n}} \egl 1 + O_P\left(\frac{k^{-1/2}}{\log (n)}\right).
$$
Consequently, 
the first term of (\ref{eqtmp1}) 
converges in distribution to ${\mathcal N}(0,\theta^2\sigma^2(W)/\mu^2(W))$. 
We also have that the second term of~(\ref{eqtmp1}) converges to $\lambda$
in probability and the result is proved.
\CQFD

\paragraph{\bf Proof of Corollary \ref{normalitezipf}.}

First remark that (\ref{defthetazipf}) can be rewritten as
$$
\hat{\theta}_n^{Z} =
{\displaystyle\somi \ain^Z (\log(X_{n-i+1,n})  - \log(X_{n-k+1,n})) }
\left/{\displaystyle\somi \ain^Z (\ld{n/i} - \ld{n/k})} \right.,
$$
where
\begin{eqnarray*}
\ain^Z &=& \log(n/k) \left ( \ld{n/i} - \zeta_n \right )\\
&=&\log(n/k) \left ( \log \left (1+\frac{\log(k/i)}{\log(n/k)} \right )-\tau_n \right )\\
&=& \log(k/i) + O \left ( \frac{\log^2(k)}{\log(n)}\right ) -\log(n/k)\tau_n,\\
&=& \log(k/i)-1 + O \left ( \frac{\log^2(k)}{\log(n)}\right ) + O\left(\frac{\log(k)}{k}\right),
\end{eqnarray*}
uniformly on $i=1,\ldots,k$ with Lemma~\ref{lemtaun}.
Therefore, we have $\ain^Z=W(i/k)+\ein$ with $W(x)=-(\log(x)+1)$ and $\ein=O(\log^2(k)/\log(n))+O(\log(k)/k)$,
uniformly on $i=1,\ldots,k$. Then, it is easy to check that $W$ satisfies conditions \At~and \Aq~and that
$ \mu(W)= 1 $ and $ \sigma^2(W) =  2$.
\CQFD

\begin{table}[h]
\begin{center}
$
\begin{array}{|c|c|c|c|c|}
\hline
                &  1-F(x) & \theta & b(x) & \rho \\
\hline
\hline
\mbox{Weibull }     & \exp(-(x/\lambda)^\alpha )& 1/\alpha & 0 & -\infty \\
{\mathcal W}(\alpha,\lambda)&&&&\\
\hline
\mbox{Extended Weibull } & r(x)\exp(-\beta
x^\tau)        & 1/\tau & -\frac{\gamma}{\tau^2}
\frac{\log x}{x} &  -1\\
{\mathcal EW}(\tau,\beta,\gamma)& r\in{\mathcal R}_\gamma  &&&\\
\hline
\mbox{Gaussian }& \frac{1}{(2\pi\sigma^2)^{1/2}}\int_x^\infty \exp\left(-\frac{(t-\mu)^2}{2\sigma^2}\right)dt    & 1/2   & \frac{1}{4} \frac{\log x}{x} & -1 \\
{\mathcal N}(\mu,\sigma^2)&&&&\\
\hline
\mbox{Gamma }& \frac{\beta^\alpha}{\Gamma(\alpha)} \int_x^\infty
 t^{\alpha-1} \exp(-\beta t) dt& 1      & (1-\alpha) 
\frac{\log x}{x} & -1 \\
\Gamma(\beta,\alpha)&&&&\\
\hline
\mbox{Benktander II }& x^{\tau-1} \exp\left(-\frac{\alpha}{\tau} x^\tau\right)& 1/\tau      & \frac{(1-\tau)}{\tau^2} \frac{\log x}{x} & -1 \\
{\mathcal B}(\alpha,\tau)&&&&\\
\hline
\mbox{Logistic }& \frac{2}{1+\exp x}   & 1     &
- \frac{\log 2}{x} & -1 \\
{\mathcal L}&&&&\\
\hline
\mbox{Extreme Value }& 1-\exp(-\exp(\mu-x))  & 1     & -\frac{\mu}{x}
 &  -1\\
{\mathcal EVD}(\mu)&&&&\\
\hline
\end{array}
$
\end{center}
\caption{Some Weibull tail-distributions}
\label{tabex}
\end{table}

\clearpage



\begin{thebibliography}{99}

\bibitem{order} Arnold, B.C, Balakrishnan, N., Nagaraja H.N., (1992),
{\em A First course in order statistics}, Wiley and sons. 

\bibitem{BAIN} Bain, L., Engelhardt, M., (1991),
{\em Statistical analysis of reliability and life-testing models. Theory and methods. 2nd ed.},
 Statistics: Textbooks and Monographs Series, {\bf 115}, Marcel Dekker, Inc.

\bibitem{BEIR} Beirlant, J., Dierckx, G., Goegebeur, Y., Matthys, G., (1999),
Tail index estimation and an exponential regression model,
{\em Extremes}, {\bf 2}, 177--200.

\bibitem{BeirBro} Beirlant J., Broniatowski, M., Teugels, J.L., Vynckier, P. (1995),
The mean residual life function at great age: Applications to tail estimation,
{\em Journal of Statistical Planning and Inference}, {\bf 45}, 21--48. 

\bibitem{BeirTeug} Beirlant, J., Teugels, J.L., (1992),
Modeling large claims in non-life insurance.
{\em Insurance: Math. Econom.}, {\bf 11}, 17--29.

 
\bibitem{BTV} Beirlant, J., Teugels, J., Vynckier, P., (1996),
{\em {Practical analysis of extreme values}},
Leuven university press, Leuven.


\bibitem{Berred} Berred, M., (1991),
Record values and the estimation of the Weibull tail-coefficient.
{\em  Comptes-Rendus de l'Acad\'emie des Sciences}
{\bf t. 312}, S\'erie I, 943--946.

\bibitem{BING} Bingham, N.H., Goldie, C.M., Teugels, J.L., (1987), {\em Regular variation}, Encyclopedia of Mathematics and its Applications, {\bf 27}, Cambridge University Press.

\bibitem{BRO} Broniatowski, M., (1993),
On the estimation of the Weibull tail coefficient,
{\em Journal of Statistical Planning and Inference}, {\bf 35}, 349--366. 

\bibitem{CDM} Cs\"org\"o, S., Deheuvels, P., Mason, D., (1985),
Kernel estimates of the tail index of a distribution,
{\em The Annals of Statistics}, {\bf 13}, 1050--1077.

\bibitem{EMBR} Embrechts, P., Kl\"uppelberg, C., Mikosch, T., (1997), {\em Modelling extremal events}, Springer.

\bibitem{Galambos}  Galambos, J., (1987),
 {\em The Asymptotic theory of extreme order statistics},
{R.E. Krieger publishing company}.

\bibitem{Revstat} Gardes, L., Girard, S., (2006),
Comparison of {W}eibull tail-coefficient estimators,
{\em REVSTAT - Statistical Journal},
{\bf{4}}, 373--188.

\bibitem{Nous} Gardes, L., Girard, S., (2008),
Estimation of the Weibull tail-coefficient with linear combination of upper order statistics,
{\em Journal of Statistical Planning and Inference},
{\bf{138}}, 1416--1427.


\bibitem{GEL} Geluk, J.L., de Haan, L., (1987),
Regular Variation, Extensions and Tauberian Theorems.
{\em Math Centre Tracts}, {\bf 40}, Centre for Mathematics and Computer Science, Amsterdam.

\bibitem{thetaWT} Girard, S.  (2004),
A Hill type estimate of the Weibull tail-coefficient, 
{\em Communication in Statistics - Theory and Methods}, {\bf 33},
205--234.

\bibitem{Hausler} H\"ausler, E., Teugels, J.L. (1985),
On asymptotic normality of Hill's estimator for the exponent of regular
variation, {\em The Annals of Statistics}, {\bf 13}, 743--756.

\bibitem{Hill} Hill, B.M., (1975),
A simple general approach to inference about the tail of a distribution,
{\em The Annals of Statistics}, {\bf 3}, 1163--1174.

\bibitem{KLUP} Kl\"uppelberg, C., Villase\~nor, J.A., (1993),
Estimation of distribution tails - a semiparametric approach,
{\em Deutschen Gesellschaft f\"ur Versicherungsmathematik},
{\bf XXI(2)}, 213--235.

\bibitem{krares96} Kratz, M., Resnick, S., (1996),
The QQ-estimator and heavy tails,
{\em Stochastic Models}, {\bf{12}}, 699--724.

\bibitem{Mason2} Mason, D.M., (1981), 
Asymptotic normality of linear combinations of order
statistics with a smooth score function,
{\em Annals of Statistics}, {\bf 9}(4), 899--908.

\bibitem{schste96} Schultze, J., Steinebach, J., (1996),
On least squares estimates of an exponential tail coefficient,
{\em Statistics and Decisions}, {\bf 14}, 353--372.

\bibitem{Viharos} Viharos, L. (1999),
Weighted least-squares estimators of tail indices,
{\em Probability and Mathematical Statistics}, {\bf 19}, 249--265.
\end{thebibliography}
\end{document}